\documentclass[11pt, reqno]{amsart}
\usepackage{amsmath,amsfonts, mathtools,amssymb,latexsym,amssymb,amsthm,adjustbox,mathrsfs,amsbsy, bm,tikz-cd,indentfirst, makecell,graphicx, verbatim, calc, enumerate,xy}
\usepackage[top=1in,bottom=1.3in,left=1in,right=1in]{geometry}
\usepackage[colorlinks=true,linkcolor=blue]{hyperref}
\usepackage{listings}
\usepackage{xcolor}
\theoremstyle{plain}
\newtheorem{theorem}{Theorem}[section]

\newtheorem{conjecture}[theorem]{Conjecture}

\newtheorem{proposition}[theorem]{Proposition}

\newtheorem{point}[theorem]{}

\newcommand{\m}{\mathfrak{m}}

\newcommand{\n}{\mathfrak{n}}

\newcommand{\N}{\mathbb{N}}

\usepackage{graphicx}
\usepackage{subfigure}
\usepackage{epsfig}
\usepackage{epstopdf}
\usepackage{float}
\theoremstyle{plain}

\usepackage{xcolor}
\usepackage{titlesec}

\titleformat{name=\section}{}{\thetitle.}{0.5em}{\centering\scshape\large}

\begin{document}
\title{\large \textbf{Hilbert-Kunz Multiplicity and Smirnov's conjecture}}
% \author{Shreedevi Masuti}
% \email{shreedevi@iitdh.ac.in}
% \address{ Indian Institute of Technology Dharwad, Dharwad 580011, Karnataka, India}

\author{Samarendra Sahoo}
\email{samarendra.s.math@gmail.com, ra.samarendra@iitdh.ac.in}
\address{ Indian Institute of Technology Dharwad, Dharwad, Karnataka, India, 580011}
%\date{\today}


\begin{abstract}
 In this paper, we study Smirnov's conjecture: Let $(R,\mathfrak{m})$ be a Noetherian local ring of characteristic $p>0$ and dimension $d$. Let $I$ be an $\mathfrak{m}$-primary ideal. If $R$ is weakly $F$-regular, analytically unramified and Cohen-Macaulay, then $\displaystyle L_1(I)=e(I)-e_{\mathrm{HK}}(I)$ if and only if $I$ is stable. Here, $e(I)$ and $e_{\mathrm{HK}}(I)$ denote the Hilbert-Samuel multiplicity and the Hilbert-Kunz multiplicity of $R$ with respect to $I$, respectively. Moreover, $\displaystyle L_1(I)=\lim_{q\to\infty}
\frac{e_1(I^{[q]})}{q^d},$ where $q=p^e$ and $e_1(I^{[q]})$ denotes the first Hilbert coefficient of $R$ with respect to the $q$th Frobenius power of $I$.

We prove the conjecture when $R$ is either excellent and strongly $F$-regular, or excellent of dimension $3$.

\end{abstract}

\maketitle

\section{Introduction}

Let $(R,\mathfrak{m})$ be a Noetherian local ring of prime characteristic $p>0$ and dimension $d$. Let $\ell_R(M)$ denote the length of an $R$-module $M.$ For an $\mathfrak{m}$-primary ideal $I$, the Hilbert-Kunz multiplicity is defined by
$$e_{\mathrm{HK}}(I)=\lim_{q\to\infty}\frac{\ell_R(R/I^{[q]})}{q^d},$$
where $q=p^e$ and $I^{[q]}=(x^q\mid x\in I)$
denotes the $q$-th Frobenius power of $I$. The Hilbert-Samuel function of $R$ with respect to $I$ is defined as the function $H^I_R(n)= \ell_R(R/I^{n})$. It is well-known that for large values of $n$, this function can be expressed as a polynomial $P_R(n)$ of degree equal to $\dim(R)$, known as the Hilbert-Samuel polynomial of $R$  with respect to $I$. The polynomial can be written in the form $$P_R(n)=e_0(I)\binom{n+d-1}{d}-e_1(I)\binom{n+d-2}{d-1}+\cdots+(-1)^d e_d(I),$$ where $e_i(I)$ is the $i$th Hilbert coefficient of $R$ with respect to $I$. We write $e(I)$ instead of $e_0(I)$. In particular, for fixed $q=p^e$, $e_i(I^{[q]})$ is the $i$th Hilbert coefficient of $R$ with respect to $I^{[q]}$. Moreover, $e(I^{[q]})=e(I^q)=q^d e(I).$ From the expression of $P_R(n)$, we can conclude that $\displaystyle e(I)=\lim_{n\to\infty}\frac{\ell_R(R/I^{n})}{n^d}.$

The two invariants $e(I)$ and $e_{\mathrm HK}(I)$ exhibit several common properties. For example, if $R$ is a formally unmixed ring, then $e(\mathfrak m)=1$ is equivalent to $R$ being regular and also $e_{\mathrm{HK}}(\mathfrak m)=1$ equivalent to $R$ being regular. The first equivalence %between $e(\mathfrak m)=1$ and the regularity of $R$ 
is due to a classical result established by P. Samuel in the equal-characteristic case \cite{samuel} and by M. Nagata in the general case \cite{nagata}. The second equivalence %between $e_{\mathrm{HK}}(\mathfrak m)=1$ and the regularity of $R$
was established more recently by K. I. Watanabe and K. I. Yoshida \cite[Theorem 1.5]{yoshida}.

Since $e(I^{[q]})=e(I^q)=q^d e(I)$ for all $q=p^e$, we have $\displaystyle \lim_{q\to\infty} \frac{e(I^{[q]})}{q^d}=e(I).$  This naturally leads to the study of the asymptotic behavior of the Hilbert coefficients $\displaystyle \frac{e_i(I^{[q]})}{q^d}
\text{ as }q\to\infty.$ 
We set, whenever the limit exists, 
$$L_i(I) =\lim_{q\to\infty}\frac{e_i(I^{[q]})}{q^d}.$$
The existence of these limits and their relationship with the Hilbert-Kunz multiplicities of ordinary powers of $I$ have been investigated in recent work. V. Trivedi first established the existence of $L_1(I)$ for standard graded algebras and equigenerated homogeneous ideals \cite{trivedi}, while I. Smirnov later established the existence of $L_1(I)$ for every $\m$-primary ideal in a local ring \cite{smirnov} and formulated the following conjecture.

\begin{conjecture}\label{conj}
    Let $(R,\mathfrak{m})$ be a Noetherian local ring of characteristic $p>0$ and dimension $d$. Let $I$ be an $\mathfrak{m}$-primary ideal in $R$. If $R$ is weakly $F$-regular, analytically unramified and Cohen-Macaulay, then $\displaystyle L_1(I)=e(I)-e_{\mathrm{HK}}(I)$ if and only if $I$ is stable, i.e., $I^2=JI$ for a minimal reduction $J$ of $I$.
\end{conjecture}

Note that the conjecture was originally stated without the weakly $F$-regular assumption, but it is known to be false without this hypothesis (see \cite{goel}). It is also known that if $I$ is stable, then $\displaystyle L_1(I)=e(I)-e_{\mathrm{HK}}(I)$, see \cite[Theorem 1.8(2)]{watanabe}. De Stefani, Masuti, Rossi, and Verma \cite{shreedevi} proved the conjecture for dimension $1$ and $2$-dimensional excellent local rings (see \cite[Corollary 2.2, Theorem 4.9]{shreedevi}). In this paper, we establish an affirmative answer to the conjecture \ref{conj} for $d\geq 3$ under certain mild conditions. We prove

\begin{theorem}[Theorem \ref{conjecture-dim3} and \ref{conjecture-dimd}]\label{1}
    Assume the same hypothesis as in Conjecture \ref{conj}. Then the conjecture holds if either of the following conditions is satisfied.
    \begin{enumerate}
        \item  $R$ is an excellent and strongly $F$-regular ring.
         \item $R$ is an excellent ring with $d\leq 3$.
    \end{enumerate}
\end{theorem}

\noindent
\textbf{Organization of the paper:} This paper is structured into three sections. Section 2 introduces the required preliminaries. In Section 3, we first establish a lower bound for $L_1(I)-e(I)+e_{\mathrm HK}(I)$ and then give a proof of the main result (Theorem \ref{1}) using the bound.

\section{Preliminaries}
Throughout this paper, all rings are Noetherian commutative with unity and modules are finitely generated. In this section, we discuss some definitions and preliminary results that we need.

We first recall the definition of $F$-signature.
\begin{point} 
\normalfont
Let $(R,\mathfrak{m},k)$ be a Noetherian local ring of prime
characteristic $p>0$, where $k=R/\mathfrak{m}$. For each $e\geq 1$, we denote by $F_*^e(R)$ the additive group of $R$ equipped with the $R$-module structure $r\cdot F_*^e (x)=F_*^e(r^{p^e}x)
\text{ for all }r,x\in R.$
The ring $R$ is called \emph{$F$-finite} if $F_*^e(R)$ is a finitely generated $R$-module for all $e\geq 1$. We write $F_*(-)$, when $e=1$. In this case, set
\[
\alpha=\log_p\bigl(\dim_k F_*(k)\bigr)
      =\log_p[k:k^p].
\]

Let $(R,\mathfrak{m},k)$ be an $F$-finite local ring of dimension $d$.
For every $e\geq 1$, write
\[
F_*^e(R)\cong R^{\oplus a_e}\oplus M_e,
\]
where $M_e$ has no nonzero free direct summand. The integer
$a_e=a_e(R)$ is called the \emph{$e$-th Frobenius splitting number}
of $R$; it is the largest rank of a free $R$-module occurring as a
direct summand of $F_*^e(R)$. The \emph{$F$-signature} of $R$ is defined as
$$s(R)=\lim_{e\to\infty}
\frac{a_e}{p^{e(d+\alpha)}}.$$
The existence of this limit was proved in full generality by
K.~Tucker \cite{Tucker12}. %Hence, by \cite[Theorem 0.2]{AL03}, $s(R)>0$ iff $R$ is strongly $F$-regular. 
\end{point}
We now list some properties that are preserved under faithfully flat extension and useful for our main results.
\begin{point}\label{flat}
    \normalfont
 Let $\phi:(R,\mathfrak{m})\to (R',\mathfrak{m}')$ be a local ring homomorphism. Assume $R'$ is a faithfully flat $R$-algebra with $\mathfrak{m}R'=\mathfrak{m}'$. Set $I'=IR'$ for an $\m$-primary ideal $I$. Let $M$ be an $R$-module. Set $M'=M\otimes_R R'$. In these cases, it can be shown that
\begin{enumerate}
    \item $\ell_R(M)=\ell_{R'}(M')$. In particular, $\ell_R(R/I^n)=\ell_{R'}(R'/I'^n)$ for all $n\geq 1$. Thus, $e_i(I)=e_i(I')$ for all $i\geq 0.$
    %\item $H^I_M(n)(=\ell_{R}(I^nM/I^{n+1}M))=H_{M'}(n)(=\ell_{R'}({I'}^nM'/{I'}^{n+1}M'))$ for all $n\geq 0$. In particular, $e_i(I)=e_i(I')$ for all $i\geq 0.$
    \item dim$_{R}M=$ dim$_{R'}M'$ and $\text{ grade}_R (K,M)= \text{ grade}_{R'} (K',M')$ for any ideal $K$ of $R$.
  \item For $R'=\widehat{R}$(the completion of $A$ with respect to the maximal ideal) with characteristic $p>0$, the following properties hold.
  \begin{enumerate}
      \item[{(a)}] $I^{[q]}R'=I'^{[q]}$ for all $q=p^e$, $e\in \N$.
      \item[{(b)}] \cite[Corollary~7.28]{HH94} $R$ is excellent weakly $F$-regular iff $\widehat{R}$ is excellent weakly $F$-regular.
      \item[{(c)}] \cite[Theorem 2.23]{Yao} $R$ is strongly $F$-regular iff $\widehat{R}$ is strongly $F$-regular.
  \end{enumerate}
\end{enumerate}

\end{point}

Next, we list some cases where weakly $F$-regular implies strongly $F$-regular. Note that a local weakly $F$-regular ring is a domain. 

\begin{proposition}\label{weakstrong}
    Let $(R,\mathfrak{m})$ be a Noetherian local ring of characteristic $p>0$. Let $R$ be a weakly $F$-regular, $F$-finite local ring of dimension $d\leq 3$ and  a homomorphic image of a Gorenstein ring. Then $R$ is strongly $F$-regular. 
\end{proposition}
\begin{proof}
 For $d=2\text{ and }3$, the result now follows from \cite[pg-736, 742]{Wil95}. Let $d=1$. By \cite[Theorem 3.4(c)]{HH94}, $R$ is normal and hence regular. By \cite[Theorem 5.5(c)]{HH94}, the result follows. 
\end{proof}

\section{main results}
In this section, we first establish a lower bound for $L_1(I)-e(I)+e_{\mathrm HK}(I)$ in terms of the $F$-signature and Hilbert coefficients of $R$ with respect to $I$. We then use this bound, together with the $\Gamma$ construction (see \cite{HH94}), to prove Theorem \ref{1}.

We begin by establishing a lower bound for $L_1(I)-e(I)+e_{\mathrm HK}(I)$.% in terms of the $F$-signature and $e_1(I)-e(I)+\ell_R(R/I)$.

\begin{theorem}\label{Fsignature}
Let $(R,\mathfrak m,k)$ be an $F$-finite Cohen-Macaulay local
ring of characteristic $p>0$ and dimension $d$, and let
$I\subseteq R$ be an $\mathfrak m$-primary ideal. Then
\[
L_1(I)-e(I)+e_{\mathrm{HK}}(I)
\geq
s(R)\bigl(e_1(I)-e(I)+\ell_R(R/I)\bigr).
\]
\end{theorem}

\begin{proof}
Let $q=p^e$ and set $\alpha=\log_p [k:k^p].$ For all $n\geq 1$, we have $I^nF_*^e(R) =F_*^e\bigl((I^{[q]})^n\bigr).$ Since $F^e_*(-)$ is an exact functor of $R$-modules, we have
$$\frac{F_*^e(R)}{I^nF_*^e(R)}
\cong
F_*^e\left(\frac{R}{(I^{[q]})^n}\right).$$ We also have, for a finite length $R$-module $N$, $\ell_R(F_*^eN)=q^\alpha\ell_R(N)$ (see \cite[2.1]{Tucker12}). It follows that
$$\ell_R\left(\frac{F_*^e(R)}{I^nF_*^e(R)}\right)
=q^\alpha \ell_R\left(\frac{R}{(I^{[q]})^n}\right)$$
for all $n\geq1$. Comparing the Hilbert polynomials, we obtain
$e_i(I,F_*^e(R))=q^\alpha e_i(I^{[q]})$ for all $0\leq i\leq d$. In particular, we get
\begin{equation} \label{eq}
    e(I,F_*^e(R))=q^\alpha e(I^{[q]})=q^{d+\alpha}e(I) \text{ and } e_1(I,F_*^e(R))=q^\alpha e_1(I^{[q]}).
\end{equation}
For a maximal Cohen-Macaulay $R$-module $M$, set
$$\delta_I(M)=e_1(I,M)-e(I,M)+\ell_R(M/IM).$$
By the module version of Northcott's inequality (see \cite{HCCMM}, Proposition 12), we have $\delta_I(M)\geq0$. Since $R$ is Cohen-Macaulay, $F_*^e(R)$ is a maximal Cohen-Macaulay $R$-module. Indeed, if $x_1,\ldots,x_d$ is an $R$-regular sequence, then $x_1^q,\ldots,x_d^q$ is also an $R$-regular sequence, and
multiplication by $x_i$ on $F_*^e(R)$ corresponds to multiplication by $x_i^q$ on the underlying copy of $R$. For all $e\geq 1$, write $$F_*^e(R) \cong R^{\oplus a_e}\oplus M_e,$$
where $M_e$ has no free direct summands. Since $F_*^e(R)$ is maximal Cohen-Macaulay, so is $M_e$. Since the Hilbert coefficients and lengths are additive on direct sums, we have
$$\delta_I(F_*^e(R))=a_e\delta_I(R)+\delta_I(M_e)\geq a_e\delta_I(R).$$
On the other hand, from \eqref{eq},
$\displaystyle \delta_I(F_*^e(R)) =q^\alpha(e_1(I^{[q]})-q^de(I)+\ell_R(R/I^{[q]}))\geq a_e\delta_I(R)).$ Dividing by $q^{d+\alpha}$ yields
$$\frac{e_1(I^{[q]})}{q^d}-e(I)+\frac{\ell_R(R/I^{[q]})}{q^d}
\geq\frac{a_e}{q^{d+\alpha}}(e_1(I)-e(I)+\ell_R(R/I)).$$
Taking the limit as $q\to\infty$, we obtain
$$L_1(I)-e(I)+e_{\mathrm{HK}}(I) \geq s(R)(e_1(I)-e(I)+\ell_R(R/I)).$$
\end{proof}

Note that if $s(R)\neq 0$ (for example, $R$ is strongly $F$-regular, see \cite[Theorem 0.2]{AL03}) and $L_1(I)=e(I)-e_{\mathrm{HK}}(I)$ in Theorem \ref{Fsignature}, then $e_1(I)=e(I)+\ell_R(R/I)$ which is equivalent to $I$ being stable (see \cite[Theorem 16]{HCCMM}). We now prove the Conjecture \ref{conj} when $R$ is an excellent and strongly $F$-regular ring. Since stabile ideal $I$ implies $L_1(I)=e(I)-e_{\mathrm{HK}}(I)$ in general (see \cite[Theorem 1.8(2)]{watanabe}), it remains only to establish the converse.

\begin{theorem}\label{conjecture-dimd}
Assume the same hypothesis as in Conjecture \ref{conj} and let $R$ be an excellent and strongly $F$-regular ring. Then
$L_1(I)=e(I)-e_{\mathrm{HK}}(I)$ if and only if $I$ is stable.
\end{theorem}
\begin{proof}
Assume that $L_1(I)=e(I)-e_{\mathrm{HK}}(I)$. By \ref{flat}, we may assume that $R$ is complete. We first reduce the case to $F$-finite ring by the well known $\Gamma$-construction and then use Theorem \ref{Fsignature}. By Cohen's structure theorem, we have $R\cong T/\mathfrak a,$ where $T=k[[x_1,\ldots,x_n]]$ and $k$ is a coefficient field of $R.$ 

Let $\Lambda$ be a $p$-basis of $k$ (see \cite[3.1]{strong}). Choose a sufficiently small cofinite subset $\Gamma\subseteq\Lambda$ and set $S=R^\Gamma.$ By \cite[6.11, 6.13, 6.6(d)]{HH94}, the map $R\to S$ is faithfully flat and purely inseparable, and $S$ is $F$-finite local excellent domain with maximal ideal $\n=\m S$. Now by \cite[Theorem 3.4(ii)]{strong}, $S$ is strongly $F$-regular. 

Since the map $(R,\m)\to (S,\n)$ is faithfully flat with $\m S=n$, by \ref{flat} and the above argument, we may assume that $R$ is $F$-finite. Therefore, by \cite[Theorem 0.2]{AL03},  $s(R)>0$ and by Theorem \ref{Fsignature}, we have $e_1(I)=e(I)-\ell_R(R/I).$ It follow that $I$ is a stable ideal, see \cite[Theorem 16]{HCCMM}.
\end{proof}

As a consequence of Theorem \ref{conjecture-dimd}, we prove the Conjecture \ref{conj} when $R$ is an excellent local ring of dimension at most $3$. Note that De Stefani, Masuti, Rossi, and Verma \cite{shreedevi} proved the conjecture for dimension $\leq2$. We now give a proof for dimension at most $3$, which is different from their proof.

\begin{theorem}\label{conjecture-dim3}
Assume the same hypothesis as in Conjecture \ref{conj} and let $R$ be an excellent ring with $d\leq 3$. Then
$L_1(I)=e(I)-e_{\mathrm{HK}}(I)$ if and only if $I$ is stable.
\end{theorem}

\begin{proof}

Assume that $L_1(I)=e(I)-e_{\mathrm{HK}}(I).$ By \ref{flat}, we may assume that $R$ is complete. We first reduce the case to $F$-finite ring by the well-known $\Gamma$-construction and then use Theorem \ref{Fsignature}. By Cohen's structure theorem, we have $R\cong T/\mathfrak a,$ where $T=k[[x_1,\ldots,x_n]]$ and $k$ is a coefficient field of $R.$ 

Let $\Lambda$ be a $p$-basis of $k$ (see \cite[3.1]{strong}). Choose a sufficiently small cofinite subset $\Gamma\subseteq\Lambda$ and set $S=R^\Gamma.$ By \cite[6.11, 6.13]{HH94}, the map $R\to S$ is faithfully flat and purely inseparable, and $S$ is $F$-finite local excellent domain with maximal ideal $\n=\m S$. In particular, $\m S$ is radical.
Therefore, \cite[Theorem~6.17(b)]{HH94} shows that
$S$ is weakly $F$-regular. By \cite[Lemma~6.6(a,b)]{HH94}, we also have $S\cong T^{\Gamma}/\mathfrak a T^{\Gamma}$, where $T^{\Gamma}$ is a regular local ring. Thus, $S$ is a homomorphic image of a Gorenstein local ring. Hence, by Proposition \ref{weakstrong}, $S$ is strongly $F$-regular.

Since the map $(R,\m)\to (S,\n)$ is faithfully flat with $\m S=n$, by \ref{flat} and the above argument, we may assume $R$ is $F$-finite and strongly $F$-regular. The result now follows from Theorem \ref{conjecture-dimd}. 
\end{proof}

\noindent
\textbf{Acknowledgment:} The author thanks Dr. Shreedevi Masuti for suggesting the problem and for her valuable comments and discussions. The author also thanks ChatGPT (version 5.6, Sol) for its assistance in working out the proof of Theorem \ref{Fsignature}.

\providecommand{\bysame}{\leavevmode\hbox to3em{\hrulefill}\thinspace}
\providecommand{\MR}{\relax\ifhmode\unskip\space\fi MR }
 %\MRhref is called by the amsart/book/proc definition of \MR.
\providecommand{\MRhref}[2]{
  \href{http://www.ams.org/mathscinet-getitem?mr=#1}{#2}
}

\end{document}